\documentclass[11pt, letterpaper]{article}

\usepackage{amsmath}
\usepackage[export]{adjustbox}
\usepackage[margin=1in]{geometry}

\usepackage{overpic}

\newcommand{\bmat}[1]{\begin{bmatrix}#1\end{bmatrix}}
\newcommand{\prox}{\operatorname{prox}}

\newcommand{\jhkcircled}[1]{\raisebox{.5pt}{\textcircled{\raisebox{-.9pt}{#1}}}}
\newcommand{\jhkhl}[1]{#1}
\newcommand\blfootnote[1]{%
  \begingroup
  \renewcommand\thefootnote{}\footnote{#1}%
  \addtocounter{footnote}{-1}%
  \endgroup
}
\graphicspath{{./graphics_included}}

\usepackage{amssymb}  
\usepackage{mathtools}
\usepackage{standalone}
\usepackage{tikz,pgfplots,pgfplotstable}

\title{
Powered Descent Guidance via First-Order Optimization with Expansive Projection
}

\author{
Jiwoo Choi\thanks{
Graduate student, Department of Aerospace Engineering, Inha University, South Korea. \texttt{jiwoochoi@inha.edu}
}
\and
Jong-Han Kim\thanks{
Associate professor, Department of Aerospace Engineering, Inha University, South Korea. \texttt{jonghank@inha.ac.kr}  (Corresponding author)
} 
\blfootnote{
This work was supported 
in part by the KSLV-II Enhancement Program funded by the Korea government(MSIT) (RS-2022-00164702), 
and in part by the Space Challenge Program funded by the Korea government(MSIT) (NRF-2022M1A3B8074247).}
}

\begin{document}

\maketitle

\begin{abstract}

This paper introduces a first-order method for solving optimal powered descent guidance (PDG) problems, that directly handles the nonconvex constraints associated with the maximum and minimum thrust bounds with varying mass and the pointing angle constraints on thrust vectors.
This issue has been conventionally circumvented via lossless convexification (LCvx), which lifts a nonconvex feasible set to a higher-dimensional convex set, and via linear approximation of another nonconvex feasible set defined by exponential functions.
However, this approach sometimes results in an infeasible solution when the solution obtained from the higher-dimensional space is projected back to the original space, especially when the problem involves a nonoptimal time of flight.
Additionally, the Taylor series approximation introduces an approximation error that grows with both flight time and deviation from the reference trajectory.
In this paper, we introduce a first-order approach that makes use of orthogonal projections onto nonconvex sets, allowing expansive projection (ExProj). 
\jhkhl{We show that 1) this approach produces a feasible solution with better performance even for the nonoptimal time of flight cases for which conventional techniques fail to generate achievable trajectories and 2) the proposed method compensates for the linearization error that arises from Taylor series approximation, thus generating a superior guidance solution with less fuel consumption.}
\jhkhl{We provide numerical examples featuring quantitative assessments to elucidate the effectiveness of the proposed methodology, particularly in terms of fuel consumption and flight time. Our analysis substantiates the assertion that the proposed approach affords enhanced flexibility in devising viable trajectories for a diverse array of planetary soft landing scenarios.}

\end{abstract}

\section{Introduction}
Recently, the space industry has been undergoing a radical transformation with the advent of the ``new space'' era.
In particular, SpaceX's Falcon 9 and Blue Origin's New Shepard have successfully demonstrated soft landing technologies for reusable launch vehicles, and the European Space Agency is working on the Themis program~\cite{blackmore2016autonomous, vila2019technology}.
In earlier days of planetary landing missions, classical guidance algorithms such as those used in the Apollo guidance system were developed, but they were limited to missions with various constraints. 
\jhkhl{In recent times, there has been notable progress in the development of convex optimization-based soft landing methodologies}~\cite{liu2017survey, mao2018tutorial, xie2020convex, sagliano2021onboard, malyuta2022convex, guadagnini2022model, kamath2023customized, sagliano2024six, yang2024rocket, shaffer2024implementation} \jhkhl{alongside machine learning-based techniques}~\cite{song2021feasibility, wang2023real, li2023closed, ccelik2024optimal}. \jhkhl{These advancements aim to effectively address a diverse array of practical constraints encountered in the domain.}

This paper addresses optimal powered descent guidance (PDG) problems with the goal of finding minimum-fuel soft landing trajectories while considering a variety of constraints on the state variables and control forces~\cite{lossless}.
Two notable constraints regarding PDG problems are as follows: 1) the magnitude of the thrust is bounded not only by an upper limit but also by a lower limit since the main engine cannot be turned off after ignition during the powered descent phase, and 2) the thrust vector must lie within a certain range because of maneuvering requirements or the field-of-view limits of the navigation sensors.
The lower bound constraints are obviously nonconvex, and the pointing angle constraints can also be nonconvex depending on the range.
In previous works relying on techniques such as lossless convexification, the problem has been reformulated by linearizing and relaxing these nonconvex constraints by introducing a slack variable limiting the magnitude of the thrust and the changes of variables handling mass variations~\cite{lossless, LCvx_mars, LCvx_enhancements}.

It is known that the relaxation is exact, with the slack variable matching the thrust magnitude, when the minimum fuel problem is solved for cases with the optimal time of flight, $t_f^*$. 
However, when the problem is solved for a nonoptimal time of flight, the solution obtained from the relaxation may be infeasible for the original problem. 
Conventionally, the optimal time of flight is found by applying bisection or similar search algorithms; however, doing so increases the computational burden and hence is not desirable during the powered descent phase, in which fast decision making is needed.
Moreover, the previous techniques additionally require Taylor series approximation for handling the thrust bound constraints, which are reformulated in the form of exponential functions~\cite{lossless,malyuta2022convex}.
Such approximation tends to result in suboptimality or constraint violation when the problem considers a significantly long time horizon or when the obtained optimal solution deviates notably from the reference trajectory along which the dynamics are linearized.

\jhkhl{Motivated by these shortcomings inherent in existing techniques, this paper introduces a novel first-order solution algorithm that directly handles the nonconvex sets originating from the original PDG problem. Unlike prior approaches which resort to relaxation or linear approximation techniques, the proposed algorithm handles the complexities of nonconvex sets without such compromises, thereby effectively surmounting the challenges highlighted earlier.}
It is observed that for the \emph{lossless convexification} cases (with the optimal time of flight), the proposed approach finds the same optimal solution that the previously known approach finds. Additionally, the proposed approach is observed to find very good practical solutions even for the \emph{lossy convexification} cases (with a nonoptimal time of flight), for which the previous approach fails.

To the best knowledge of the authors, this is the first results in directly handling the nonconvex constraints from the PDG problems under the convex optimization frameworks. We summarize the advantages of the proposed approach as follows:
\begin{itemize}
\item It solves PDG problems with nonconvex constraints at the computational complexity of a single convex optimization problem.
\item It generates practically feasible trajectories even when the previously known approaches fail.
\item It significantly reduces the suboptimality of the solution originating from linear approximation error.
\end{itemize}

Notably, the proofs of convergence for first-order methods such as the alternating direction method of multipliers (ADMM) or proximal gradient (PG) techniques largely rely on the fact that projection onto convex sets is \emph{nonexpansive}.
However, note that nonexpansivity of the projection operators is not a necessary and sufficient condition for algorithm convergence; in fact, it is an overly conservative condition, as there are many cases in which expansive projections can still achieve global convergence.
Our approach involves orthogonal projection onto nonconvex sets, which can be \emph{expansive}, and is thus deserving of the name expansive projection (ExProj).

The rest of this article is arranged as follows. We begin with briefly describing the mathematical formulation and the nonconvex nature of the PDG problem. Then we explain how the nonconvex constraints in PDG problem can be directly handled with first-order optimization methods, and give the algorithmic details for efficiently computing the solutions.
We finally present the computational results that shows the advantages of the proposed expansive projection approach over the existing convexification technique, and we also present an indoor flight test result that justifies the real-time implementation of the proposed approach.

\begin{table}[b]
\caption{Nomenclature}
\label{table:nomenclatures}
\begin{center}
\begin{tabular}{c l l}
\hline
\hline
$r(t)$ & Position vector \\
$\dot{r}(t)$ & Velocity vector \\
$T(t)$ & Thrust vector \\
$m(t)$ & Vehicle mass \\
$\rho_1$ & Minimum required thrust \\
$\rho_2$ & Maximum allowable thrust\\
$g$ & Gravitational acceleration vector \\
$\alpha $ & Constant describing mass consumption rate\\
$\theta_{\text{tp}} $ & Maximum thrust pointing angle from the vertical \\
\hline
\hline
\end{tabular}
\end{center}
\end{table}

\section{Powered Descent Guidance Problem}

\subsection{The PDG Problem in Nonconvex Form}

The minimum-fuel PDG problem for planetary landers or reusable launchers is defined as follows.

\smallskip\underline{\bf{Problem-PDG:}}
\begin{subequations}
\label{eq:nonconvex_problem}
\begin{align}  
 \underset{t_f, T(\cdot)}{\text{minimize}} \quad & \displaystyle \int_0^{t_f} \left\|T(t) \right\|dt \nonumber \\
\text{subject to} \quad  & \ddot{r}(t) = g + T(t)/m(t), \label{eq:dynamic}  \\
&\dot{m}(t) = -\alpha \left\|T(t) \right\|, \label{eq:mass} \\
&0 < \rho_1 \leq \left\|T(t) \right\| \leq \rho_2,  \label{eq:thrust_bound}\\
&e_1^T T(t) \geq \left\|T(t) \right\| \cos{\theta_\text{tp}}, \label{eq:thrust_direction}\\
&r(0) = r_\text{init},\ \dot{r}(0) = \dot{r}_\text{init},\ m(0) = m_\text{wet}, \label{eq:state_init}\\
& r(t_f) = \dot{r}(t_f) = 0,\ m(t_f)  \geq m_\text{dry}, \label{eq:state_final}\\
&\forall t \in [0, t_f],  \nonumber
\end{align} 
\end{subequations}
where \eqref{eq:dynamic} expresses the dynamics under the assumption that the vehicle is a point mass and \eqref{eq:mass} gives the change in mass due to thrust usage.
The constraints that follow, from \eqref{eq:thrust_bound} to \eqref{eq:state_final}, describe the limits on the thrust magnitude, the pointing angle constraint, and the initial and final constraints, respectively.
In addition to the constraints described in \jhkhl{Problem-PDG}, a wide range of additional constraints can be considered depending on the mission and vehicle types.
Note that the constraints in \eqref{eq:dynamic} and \eqref{eq:mass} are nonlinear and that the lower bound on the thrust in \eqref{eq:thrust_bound} is nonconvex.
In addition, \eqref{eq:thrust_direction} can be nonconvex depending on $\theta_\text{tp}$.
Moreover, note that $t_f$ is an optimization variable rather than being fixed {\it a priori}.

\subsection{Lossless Convexification}

In the lossless convexification (LCvx) technique described in \cite{lossless, LCvx_mars, LCvx_enhancements}, a slack variable $\Gamma(t)$ such that $\|T(t)\|\le\Gamma(t)$ is introduced to apply convex relaxation to the nonconvex constraints of \jhkhl{Problem-PDG}: the lower bound on the thrust magnitude, $\rho_1 \leq \left\| T(t) \right\|$, and the pointing angle constraint on the thrust vector, $e_1^T T(t) \geq \left\| T(t) \right\| \cos{\theta_\text{tp}}$.
Additionally, with the following changes of variables for linearizing the dynamics,
\begin{equation}
\begin{aligned}
    u(t) &= T(t)/m(t), \\
    \sigma(t) &= \Gamma(t) / m(t), \\
    z(t) &= \log m(t), \\
    z_{\text{ref}}(t) &= \log( m_\text{wet} - \alpha \rho_2 t),
\end{aligned} 
\end{equation}
the original constraint on the slack variable and the thrust bounds,
$\|T(t)\|\leq\Gamma(t)$ and $\rho_1 \leq \Gamma(t) \leq \rho_2$, respectively, are reformulated as follows.
\begin{equation}
\label{eqn:reformulated_bounds_slack}
\|u(t)\| \le \sigma(t),
\end{equation}
\begin{equation}
\label{eqn:ncvx_thrust}
\rho_1 e^{-z(t)} \leq \sigma (t) \leq \rho_2 e^{-z(t)}.
\end{equation}

Note that in \eqref{eqn:ncvx_thrust}, the inequality on the left-hand side defines a nonlinear but convex set, while the inequality on the right-hand side defines a nonconvex set. LCvx linearizes both sides around a reference trajectory $z_{\text{ref}}(t)$, yielding the following convex problem.

\smallskip\underline{\bf{Problem-LCvx:}}
\begin{equation}
\label{eqn:LCvx}
\begin{aligned}  
 \underset{t_f, \sigma(\cdot), u(\cdot)}{\text{minimize}} \quad & \displaystyle \int_0^{t_f} \sigma (t)dt \\
\text{subject to} \quad  & \ddot{r}(t) = g + u(t),  \\
&\dot{z}(t) = -\alpha\sigma (t),  \\
&\left\|u(t) \right\| \leq \sigma (t),  \\
&e_1^T u(t) \geq\sigma (t) \cos{\theta}_\text{tp},  \\
& \rho_1 e^{-z_\text{ref}(t)}\{1 - (z(t) - z_\text{ref}(t))\} \leq \sigma(t) ,\\
& \sigma(t) \leq \rho_2 e^{-z_\text{ref}(t)}\{1 - (z(t) - z_\text{ref}(t))\},  \\
&r(0) = r_\text{init},\ \dot{r}(0) = \dot{r}_\text{init},\ m(0) = m_\text{wet}, \\
& r(t_f) = \dot{r}(t_f) = 0,\ m(t_f)  \geq m_\text{dry},  \\
&\forall t \in [0, t_f].
\end{aligned} 
\end{equation}

Note that \jhkhl{Problem-LCvx} is equivalent to \jhkhl{Problem-PDG} only when the inequality in \eqref{eqn:reformulated_bounds_slack} is tight.
\begin{equation}
\label{eqn:ncvx_slack}
\|u(t)\| = \sigma(t).
\end{equation}

It has been proven that the optimal solution obtained from \jhkhl{Problem-LCvx} for the optimal time of flight of \jhkhl{Problem-PDG}, $t_f^*$, is guaranteed to satisfy $\left\|u^*(t) \right\| = \sigma^*(t), \forall t \in [0, t_f^*] $.

\section{First-Order Optimization}
\subsection{First-Order Methods}
First-order methods are a class of optimization algorithms that use the gradient or subgradient information along with proximal operations to solve optimization problems.
Compared to the classical second-order methods such as interior point methods~\cite{boyd2004convex,nemirovski2008interior,gondzio2012interior}, which require the Hessian information, first-order methods are computationally simpler and more robust, easier to implement, and more efficient in terms of memory usage; hence, they can handle problems on very large scales~\cite{beck2017, ADMM, ryu2022}.

Although first-order methods can be less accurate and may require more iterative computations than second-order methods, they are well suited for control applications in which a high-accuracy solution is not necessarily required and there are successive opportunities to update the solution in subsequent time steps~\cite{PIPG, cone_projection}.

\subsection{Convergence, Nonexpansivity, and Expansive Projection}

A general convex optimization problem with convex $f(\cdot)$ and $\mathcal{C}$ can be expressed as
\begin{equation}
	\underset{x\in\mathcal{C}}{\text{minimize}} \quad f(x)
\end{equation}
and can be solved via the ADMM, for example, by iteratively updating
\begin{equation}
\label{eqn:admm}
\begin{aligned}
	x^{k+1} &= \prox_{f/\rho}\left(z^{k}-u^k\right), \\
	z^{k+1} &= \prox_{I_\mathcal{C}}\left(x^{k+1}+u^k\right), \\
	u^{k+1} &= u^k +  x^{k+1} - z^{k+1},
\end{aligned}
\end{equation}
where $I_\mathcal{C}(\cdot)$ denotes the indicator function of the set $\mathcal{C}$ defined by
\begin{equation}
I_{\mathcal{C}}(x) = \begin{cases}
    0, & \text{if} \  x \in \mathcal{C} \\
    \infty, & \text{otherwise}
\end{cases}
\end{equation}

The global convergence of \eqref{eqn:admm} can be shown based on an understanding of the nonexpansivity of the proximal operator of convex functions~\cite{beck2017, nonexpansive, ryu2022}.
Noting that the orthogonal projection $\Pi_\mathcal{C}(\cdot)$ is the proximal operator of the indicator function $I_\mathcal{C}(\cdot)$, nonexpansivity for the convex $f(\cdot)$ and $\mathcal{C}$ implies that
\begin{equation}
\label{eqn:nonexpansivity-a}
\| \prox_f(x) - \prox_f(y) \| \le \|x-y\| 
\end{equation}
and
\begin{equation}
\label{eqn:nonexpansivity-b}
\| \Pi_\mathcal{C}(x) - \Pi_\mathcal{C}(y) \| \le \|x-y\|
\end{equation}
for all $x$ and $y$.

\begin{figure}[t]
	\centering
    \includegraphics[width=0.6\linewidth]{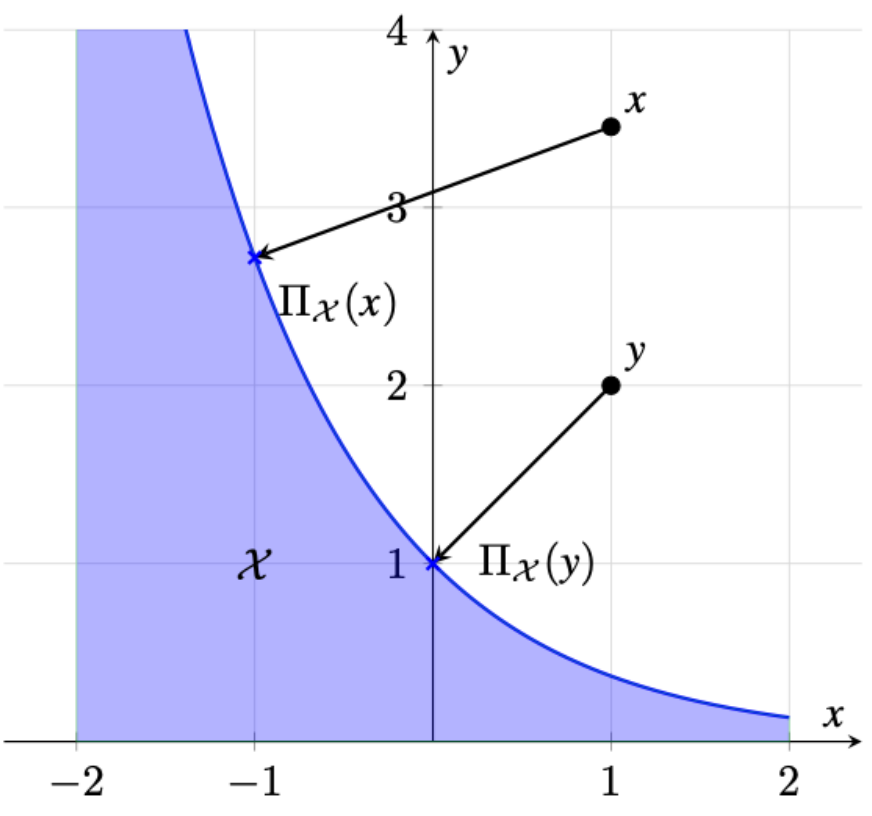}
\caption{Orthogonal projections onto a nonconvex set $\mathcal{X}$. \jhkhl{Observe that the projection onto nonconvex sets can be expansive.}}
    \label{fig:nonexpansive}
\end{figure}

On the other hand, as shown in Figure~\ref{fig:nonexpansive}, the operator $\Pi_{\mathcal{X}}(\cdot)$ for projection onto a nonconvex set $\mathcal{X}$ can be expansive, that is
\begin{equation}
    \left\| \Pi_{\mathcal{X}}(x) - \Pi_{\mathcal{X}}(y) \right\| \geq \left\| x - y \right\|,
\end{equation}
for some $x$ and $y$.

Note that projection onto nonconvex sets can sometimes increase the distance between two points, violating the nonexpansivity condition in \eqref{eqn:nonexpansivity-b}. However, this does not necessarily imply divergence of \eqref{eqn:admm}, as \eqref{eqn:nonexpansivity-b} is merely a sufficient but not necessary condition for the convergence of \eqref{eqn:admm}. This consideration can serve as a logical background for the intuitive application of expansive projections to first-order optimization.
\jhkhl{Although general results on the convergence of the first-order methods on nonconvex problems are not yet known in general, however readers can refer to}~\cite{wang2019global, boct2020proximal, zhang2020proximal, le2020inertial, themelis2020douglas, yi2021linear, danilova2022recent, hurault2022proximal, yang2023proximal, boct2023alternating} 
\jhkhl{for some recently reported analytic results on nonconvex optimization.}

\subsection{First-Order Method with Expansive Projection for PDG Problems}

The PDG problem convexified via LCvx can be readily solved using off-the-shelf solvers. However, this approach is not able to produce a practically feasible solution for some cases with $t_f$ not equal to $t_f^*$; moreover, the solution may be suboptimal due to the approximation error arising from linearizing \eqref{eqn:ncvx_thrust}.

In this paper, we do not apply convex relaxation or linear approximation. Instead, we use a first-order method with direct projection onto nonconvex sets, which can be expansive.
With the nonconvex constraints in \eqref{eqn:ncvx_thrust} and \eqref{eqn:ncvx_slack}, the problem can be defined as follows. Note that this problem is identical to the original PDG problem in \eqref{eq:nonconvex_problem} and that two nonconvex constraints are present in \eqref{eqn:socb} and \eqref{eqn:expn}.

\smallskip\underline{\bf{Problem-ExProj:}}
\begin{subequations}
\label{eqn:problem_ExProj}
\begin{align}  
 \underset{t_f, \sigma(\cdot), u(\cdot)}{\text{minimize}} \quad & -z(t_f) \nonumber \\
\text{subject to} \quad  & \ddot{r}(t) = g + u(t), \label{eqn:ExProj_dynamic}  \\
&\dot{z}(t) = -\alpha\sigma (t), \label{eqn:Exproj_mass} \\
&\left\|u(t) \right\| = \sigma (t), \label{eqn:socb}  \\
&e_1^T u(t) \geq\sigma (t) \cos\theta_\text{tp},  \label{eqn:pointing} \\
& \rho_1 e^{-z(t)} \leq \sigma(t) \leq \rho_2 e^{-z(t)}, \label{eqn:expn}\\
&r(0) = r_\text{init},\ \dot{r}(0) = \dot{r}_\text{init},\ m(0) = m_\text{wet}, \\
& r(t_f) = \dot{r}(t_f) = 0,\ m(t_f)  \geq m_\text{dry}.  \\
&\forall t \in [0, t_f]. \nonumber
\end{align}
\end{subequations}

\subsection{ADMM with Expansive Projection}
We apply the ADMM procedures with expansive projections (ExProj) onto the nonconvex sets defined by \eqref{eqn:socb} and \eqref{eqn:expn}.
The ADMM algorithm combines dual ascent and the method of multipliers to find the optimal solution by alternately updating the primal variables and the dual variables~\cite{ADMM}.

With $x=\bmat{r^T & \dot{r}^T}^T$ and the terminal constraints expressed in the soft constraint term, Problem-ExProj can be discretized into the following standard form:
\begin{equation}
\label{eqn:descretized_problem}
\begin{aligned}
\text{minimize} \quad & -z_N + \gamma \left\|x_N \right\|^2   \\
\text{subject to} \quad & P y \geq q,  \\
           & y \in \mathcal{C}_0 \cap \mathcal{C}_1 \cap {C}_2,
\end{aligned} 
\end{equation}
with
\begin{equation}
\label{eqn:define_constrained_set}
\begin{aligned}
    & \mathcal{C}_0 = \left\{y\ \vert\ Gy = b \right\}, \\
    & \mathcal{C}_1 = \cap_i \left\{(u_i, \sigma_i)\ \vert\ \left\|u_i \right\| = \sigma_i \right\},  \\
    &  \mathcal{C}_2 = \cap_i \left\{(z_i,\sigma_i)\ \vert\ 
   \rho_1 e^{-z_i} \leq \sigma_i \leq \rho_2 e^{-z_i}\right\},
\end{aligned} 
\end{equation}
where $\gamma$ is some positive weighting parameter, and the horizon size $N$ satisfies $N\Delta t = t_f$ with sampling interval $\Delta t$.
The vector $y_i$ consists of the $i$-th state variables and control vector and takes the form $y_i = \bmat{u_i^T & x_{i+1}^T & \sigma_i & z_{i+1} }^T$ for $i \in \left\{0,1,\cdots, N-1 \right\}$, and we correspondingly define the stacked variable $y = \bmat{y_0^T & y_1^T & \cdots & y_{N-1}^T}^T$.
Thus, we can encode the inequality constraints in \eqref{eqn:pointing} as $Py \geq q$ and the dynamic constraints in \eqref{eqn:ExProj_dynamic} and \eqref{eqn:Exproj_mass} as $Gy=b$.

Accordingly, the problem is further reformulated as follows:
\begin{equation}
\label{eqn:ExProj_ADMM}
\begin{aligned}
\text{minimize} \quad & -z_N + \gamma\left\|x_N \right\|^2 + I_{\mathcal{C}_0}(y) + I_{\mathcal{C}^w}(w) \\
 \text{subject to} \quad   & D^{w_1}y - w_1 = 0, \\
    & D^{w_2}y - w_2 -b^{w_2} = 0, \\
    & Py -  q - w_3 = 0, 
\end{aligned} 
\end{equation}
with
\begin{equation}
\label{eqn:set_stacked_var}
\begin{aligned}
    \mathcal{C}_3 &= \left\{w_3 \ \vert\ w_3 \geq 0 \right\}, \\
    \mathcal{C}^w &= \left\{w \ \vert\ w_1 \in \mathcal{C}_1,\ w_2 \in \mathcal{C}_2,\ w_3 \in \mathcal{C}_3\right\}, 
\end{aligned}
\end{equation}
where $w_1$ and $w_2$ are auxiliary variables of $y$ and $D^{w_1}$ and $D^{w_2}$ are selection matrices that extract $w_1$ and $w_2$ from $y$. The variable $w_3$ is a slack variable introduced for handling the inequality constraints in $Py - q \geq 0$, and we additionally define the stacked variable $w = \bmat{w_1^T & w_2^T & w_3^T}^T$.
\begin{equation}
\begin{aligned}
     D^{w_1}y &= \begin{bmatrix}
        u_0^T & \sigma_0& u_1^T& \sigma_1 &\cdots & u_{N-1}^T &\sigma_{N-1}
    \end{bmatrix}^T, \\
     D^{w_2}y &= \begin{bmatrix}
        \sigma_0& \sigma_1& z_1& \sigma_2& z_2 &\cdots&  \sigma_{N-1}& z_{N-1}
    \end{bmatrix}^T, \\
     b^{w_2} &= \begin{bmatrix}
        z_0 & 0 & \cdots & 0
    \end{bmatrix}^T.
\end{aligned}
\end{equation}

With $H$ and $h$ defined such that $y^THy+h^Ty=-z_N + \gamma\left\|x_N \right\|^2$,  
we can formulate the augmented Lagrangian for Problem \eqref{eqn:ExProj_ADMM} as
\begin{equation}
L_{\rho}(y, w, y_s)  =
\overbrace{ y^THy + h^Ty + \frac{\rho}{2} \left\|Cy - \Tilde{q} - w + y_s \right\|_2^2}^{f_g(y,w,y_s)} 
+ I_{\mathcal{C}_0}(y) + I_{\mathcal{C}^w}(w)
\end{equation}
with
\begin{equation}
w = \begin{bmatrix}
w_1\\ w_2\\ w_3
\end{bmatrix}, \, C = \begin{bmatrix}
 D^{w_1}\\ D^{w_2}\\ P
\end{bmatrix}, \, \Tilde{q} = 
\begin{bmatrix}
0\\ b^{w_2} \\ q
\end{bmatrix}, \, y_s = \begin{bmatrix}
    y_{s,1} \\ y_{s,2} \\ y_{s,3}
\end{bmatrix},
\end{equation}
from which the problem in~\eqref{eqn:ExProj_ADMM} can be solved by iteratively applying the following updates.
\begin{equation}
\begin{aligned}
   y^{j+1} =& \
    \underset{y}{\text{argmin}} \left(f_g(y,\ w^j,\ y_s^j) + I_{\mathcal{C}_0}(y) \right), \\
    w^{j+1} =& \
    \underset{w}{\text{argmin}} \left (  
    \frac{\rho}{2}
    \left\| Cy^{j+1} - \Tilde{q} -w  + y_s^j
    \right\|^2 + I_{ \mathcal{C}^w}(w) \right), \\
    y_s^{j+1} =& Cy^{j+1} - \Tilde{q} - w^{j+1}.
\end{aligned}
\end{equation}

We handle the linear constraints $Gy = b$ without introducing an additional multiplier. Instead, we make use of the KKT condition for the associated equality-constrained quadratic problem~\cite{ADMM_update_y} as follows.
\begin{equation}
\begin{matrix}
    y^{j+1} = \
    \underset{y}{\text{argmin}} \left(f_g(y,\ w^j,\ y_s^j) + I_{\mathcal{C}_0}(y) \right) \\
    \Big\Updownarrow \\
     \begin{bmatrix}
            \rho C^TC + H & G ^T \\
            G & 0 \\
            \end{bmatrix} \begin{bmatrix}
                y^{j+1} \\ *
            \end{bmatrix}  = 
           \begin{bmatrix}
    \rho C^T\left (\Tilde{q} + w^j
 - y_s^j \right ) - h \\
    b
    \end{bmatrix}.
\end{matrix}
\end{equation}

The updates for $w_1$ and $w_2$ involve a series of expansive projection operations onto $\mathcal{C}_1$ and $\mathcal{C}_2$, respectively, while the updates for $w_3$ is relatively simple.
\begin{equation}
\begin{aligned}
     w^{j+1} =& \
    \underset{w}{\text{argmin}} \left (  
    \frac{\rho}{2}
    \left\| Cy^{j+1} - \Tilde{q} -w  + y_s^j
    \right\|^2 + I_{ \mathcal{C}^w}(w) \right) \\
    = & \begin{bmatrix*}[l]
        \Pi_{\mathcal{C}_1} \left ( D^{w_1} y^{j+1} - w_1 + y^j_{s,1} \right )\\
        \Pi_{\mathcal{C}_2} \left ( D^{w_2} y^{j+1} - w_2 + y^j_{s,2} \right ) \\
        \Pi_{\mathcal{C}_3} \left ( P y^{j+1} - q + y^j_{s,3} \right )  
    \end{bmatrix*}.
\end{aligned}
\end{equation}

\subsubsection{Projection onto $\mathcal{C}_1$}

In general, projection onto the convex set defined by a second-order cone constraint of the form $\left\|u\right\| \leq \sigma$ involves projection from the exterior region of the cone~\cite{beck2017,cone_projection}.
However, in our problem, projection from both the exterior and interior regions should be considered because the constraint in \eqref{eqn:socb} represents the surface of a second-order cone. 
\begin{equation}
\Pi_{\mathcal{C}_1}(u, \sigma) = \begin{cases}
    (u, \sigma) & \text{if } \sigma = \| u \| \\
    0 & \text{if } \sigma \leq -\| u \| \\
    \left(\frac { \| u \| + \sigma}{2\| u \|} u, \frac { \| u \| + \sigma}{2} \right) & \text{otherwise}.
\end{cases} 
\end{equation}

\begin{figure}
	\label{fig:SOCB}
	\centering
	\begin{overpic}[width=0.6\linewidth]{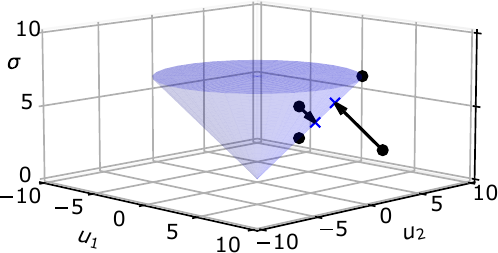}
	\put(78,21){\jhkcircled{1}}
	\put(60,17){\jhkcircled{2}}
	\put(55,32){\jhkcircled{3}}
	\end{overpic}
\caption{Orthogonal projection onto the \emph{surface of the second-order cone} defined by $\mathcal{C}_1$. Note that the set is nonconvex and that the projection from region \jhkcircled{1} (exterior region) is nonexpansive, while the projection from region \jhkcircled{3} (interior region) can be expansive.} 
\end{figure}

Note from Figure \ref{fig:SOCB} that the projection from region \jhkcircled{3} can be expansive.

\begin{figure}
	\label{fig:exp_proj}
	\centering
	\begin{overpic}[width=0.6\linewidth]{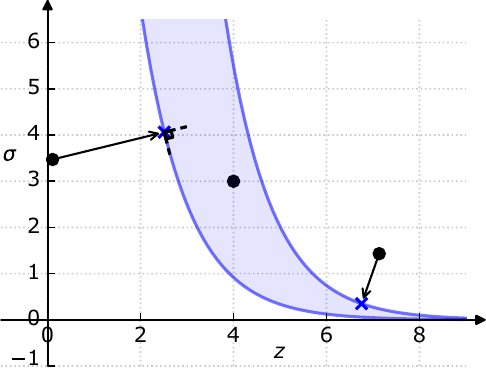}
	\put(12,38){\jhkcircled{1}}
	\put(42,32){\jhkcircled{2}}
	\put(80,24){\jhkcircled{3}}
	\end{overpic}
\caption{Orthogonal projection onto $\mathcal{C}_2$ (shaded area). Note that the set is nonconvex and that the projection from region \jhkcircled{1} \jhkhl{(below $\mathcal{C}_2$)} is nonexpansive, while the projection from region \jhkcircled{3} \jhkhl{(above $\mathcal{C}_2$)} can be expansive.} 
\end{figure}

\begin{figure}
	\centering
	\begin{overpic}[width=0.6\linewidth]{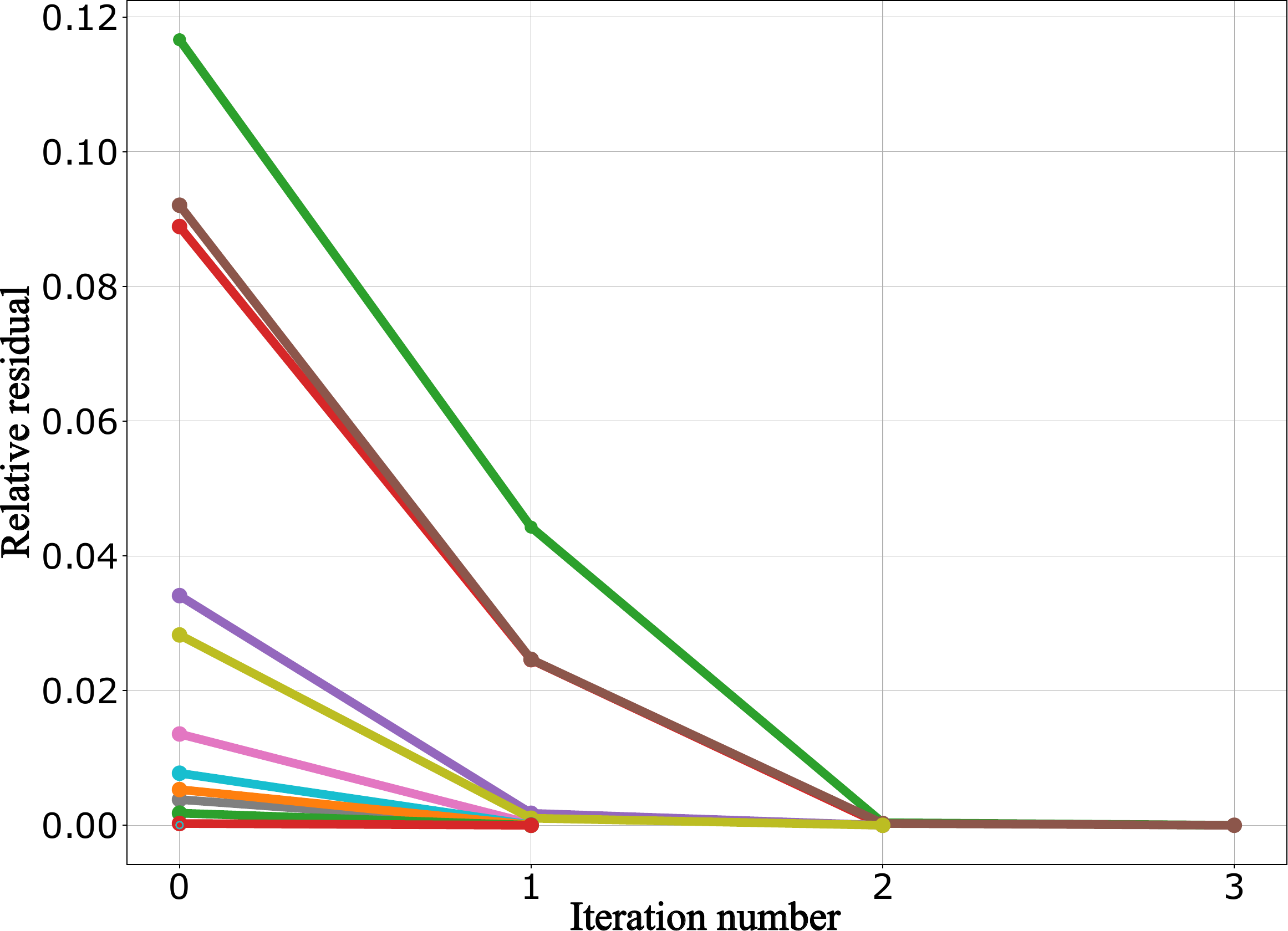}
	\end{overpic}
\caption{\jhkhl{Relative residuals attained by using the Newton-Raphson method for the instances derived from ExProj. Note that the computational process exhibits exponential convergence, typically converging within a minimal number of iterations.}}
      \label{fig:NR_iteration}
\end{figure}

\subsubsection{Projection onto $\mathcal{C}_2$}

The set $\mathcal{C}_2$ describes an area surrounded by two exponential curves.
For a point $(z,\sigma)\notin \mathcal{C}_2$, the line passing through the corresponding projected point and itself is perpendicular to the slope at the projected point on the boundary. Based on this observation, we can compute the projection as follows. 
\begin{equation}
    \Pi_{\mathcal{C}_2}(z,\sigma) = 
    \begin{cases}
    (z,\sigma) & \text{if } (z,\sigma)\in\mathcal{C}_2 \\
    (t_1,\rho_1 e^{-t_1}) & \text{if } \sigma < \rho_1 e^{-z} \\ 
    (t_2,\rho_2 e^{-t_2}) & \text{otherwise}.
\end{cases} 
\end{equation}

Here, $t_1$ and $t_2$ satisfy the followings.
\begin{equation}
\label{eqn:exponential_projection}
\begin{aligned}
&e^{t_1}(t_1-z) -\rho_1^2e^{-t_1} +\rho_1 \sigma = 0, \\
&e^{t_2}(t_2-z) -\rho_2^2e^{-t_2} +\rho_2 \sigma = 0.
\end{aligned}
\end{equation}

Note from Figure~\ref{fig:exp_proj} that the projection from region \jhkcircled{3} can be expansive.

Since the solution to \eqref{eqn:exponential_projection} is not given in closed form, we numerically solve the problem using the Newton--Raphson method, which rapidly converges within a few steps for this problem; see Figure~\ref{fig:NR_iteration}. 

\subsubsection{Projection onto $\mathcal{C}_3$}

The projection onto the set $\mathcal{C}_3$ is simply obtained by taking only the positive parts.
\begin{equation}
    \Pi_{\mathcal{C}_3} ( P y^{j+1} - q + y^j_{s,3} ) 
    = ( P y^{j+1} - q + y^j_{s,3} )_+.
\end{equation}

\section{Numerical Examples}
To verify the performance of the proposed approach (ExProj) and to make quantitative comparisons with the previously known standard technique (LCvx), we prepared three simulation cases.
The first is the case of the optimal time of flight, $t_f=t^*_f$, for which convexification with the LCvx technique is lossless and both ExProj and LCvx are expected to find the optimal solution.
The other two cases concern nonoptimal time of flight, $t_f<t^*_f$ and $t_f>t^*_f$, for which the convexification can be lossy.

The simulation scenario is based on the Mars soft landing mission presented in \cite{lossless}, and the simulation parameters are summarized in Table~\ref{table:parameters}.
The LCvx problem was solved using the {\tt ECOS}~\cite{ecos} solver via the {\tt cvxpy}~\cite{cvxpy} parser, and the ExProj problem was solved using a customized solver based on the algorithms presented in this paper.

\begin{table}[b]
\caption{Simulation Parameters} \label{table:parameters}
\begin{center}
\begin{tabular}{l c l l c}
\hline\hline
Parameter & Value & & Parameter & Value \\
\hline
$r_{\text{init}}$ (m) & (2400, 450, -330) & & $m_{\text{wet}}$ (kg) & 2000 \\
$\dot{r}_{\text{init}}$ (m/s) & (-10, -40, 10) & & $m_{\text{dry}}$ (kg) & 1700 \\
$g$ (m/s$^2$)& (-3.71, 0, 0) & & $\rho_1$ (N) & 4800\\
$\alpha$ (s/m) & 0.0005 & & $\rho_2$ (N) & 19200 \\
$\theta_{\text{tp}}$ (deg) & 90 & & & \\
\hline
\hline
\end{tabular}
\end{center}
\end{table}

\begin{figure}
	\centering
    \includegraphics[width=0.97\linewidth]{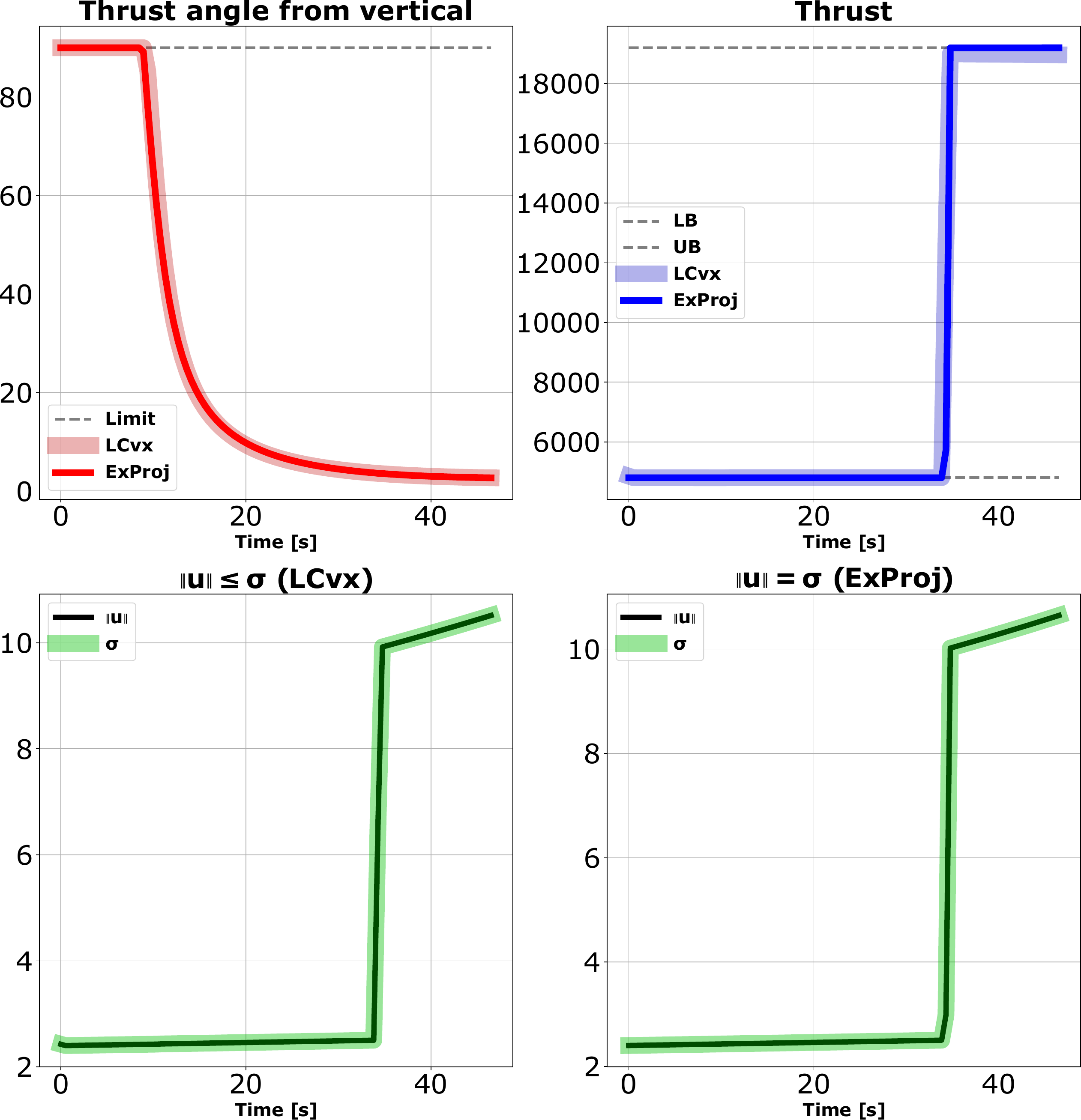}
\caption{ \jhkhl{Optimization results for the scenario of optimal time of flight case} ($t_f=t^*_f$). \jhkhl{It is observed that the convexification is lossless in this case, resulting in both LCvx and ExProj algorithms yielding nearly identical optimal solutions. However, note that in the latter phase of the thrust profile, the LCvx solution marginally underutilizes the maximum allowed thrust, rendering it slightly suboptimal.}} \label{fig:optimal_tof}
\end{figure}

\subsection{Optimal Time of Flight Case ($t_f=t^*_f$)}
For the case with the optimal time of flight, $t_f=t^*_f$, the convexification technique applied in LCvx is lossless, yielding $\|u_i\|=\sigma_i$ for all $i=0,\dots,N-1$, and both LCvx and ExProj converge to the globally optimal solution. The optimal time of flight for the given scenario is found to be $t_f^*=46.96\ \text{s}$.
The results are summarized in Figure~\ref{fig:optimal_tof} and Table~\ref{table:results}, which presents the final position and velocity as well as the fuel consumption for the mission.

We observe that the two methods produce mostly the same results, with slight differences that can be attributed to the approximation error due to the linearization of \eqref{eqn:ncvx_thrust} in LCvx.
The effect of this approximation error is also visible in the later part of the thrust profile, where the ExProj solution utilizes the full maximum thrust, while the LCvx solution uses slightly less than the maximum allowed thrust, making the LCvx solution slightly conservative and suboptimal; see Table~\ref{table:results}.
However, note that the overall effect of these differences on the achieved fuel consumption is not significant in this case.

\begin{table}
\caption{Optimization Results for $t_f = t_f^*$} \label{table:results}
\begin{center}
\begin{tabular}{l c c}
\hline
\hline
& LCvx & ExProj \\
\hline
$r$ ($10^{-5}$m)  & (0.978, -4.55, -3.93)  & (0.959, -5.13, -3.89)  \\
$\dot{r}$ ($10^{-1}$m/s)  & (-0.367, -0.412, 1.67)  & (-0.363, -0.355, 1.66) \\
$ m_{\text{wet}} - m({t_f^*})$ (kg) & 201.00  & 200.66  \\
\hline
\hline
\end{tabular}
\end{center}
\end{table}

\subsection{Nonoptimal Time of Flight Cases ($t_f<t^*_f$ or $t_f>t^*_f$)}
Cases with $t_f\neq t_f^*$ are where the convexification in the LCvx technique fails or the linearization error is significantly large.

\begin{figure}
	\centering
    \includegraphics[width=0.97\linewidth]{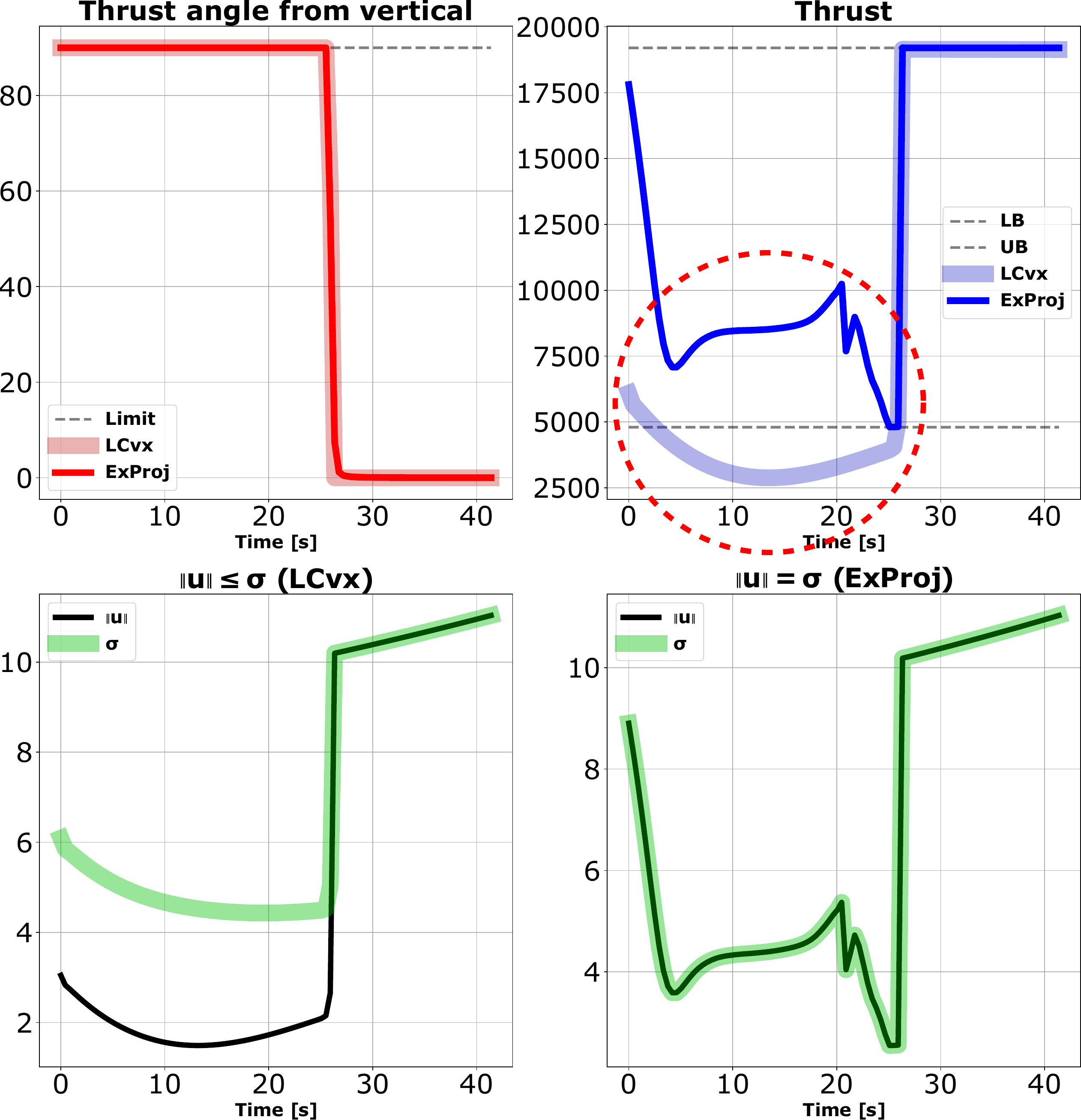}
\caption{\jhkhl{Optimization results for the scenario wherein the time of flight is strictly shorter than the optimal time} ($t_f < t_f^*$).  
\jhkhl{It is noteworthy that under these conditions, the convexification technique from the LCvx approach is no longer lossless, and the obtained solution violates the lower bound limit on the thrust (highlited in dotted red circle). Conversely, the ExProj successfully finds a feasible solution.}} \label{fig:non_optimal_tof_less}
\end{figure}

For the first case, in which the time of flight is strictly shorter than the optimal time ($t_f = 41.8\,\text{sec} < t_f^*$), we observe that the convexification due to LCvx fails with $\|u_i\|<\sigma_i$, resulting in an infeasible solution. This can be seen in Figure~\ref{fig:non_optimal_tof_less}, where the thrust profile found with LCvx is below the lower limit, whereas the ExProj solution is still feasible.

\begin{figure}
	\centering
    \includegraphics[width=0.97\linewidth]{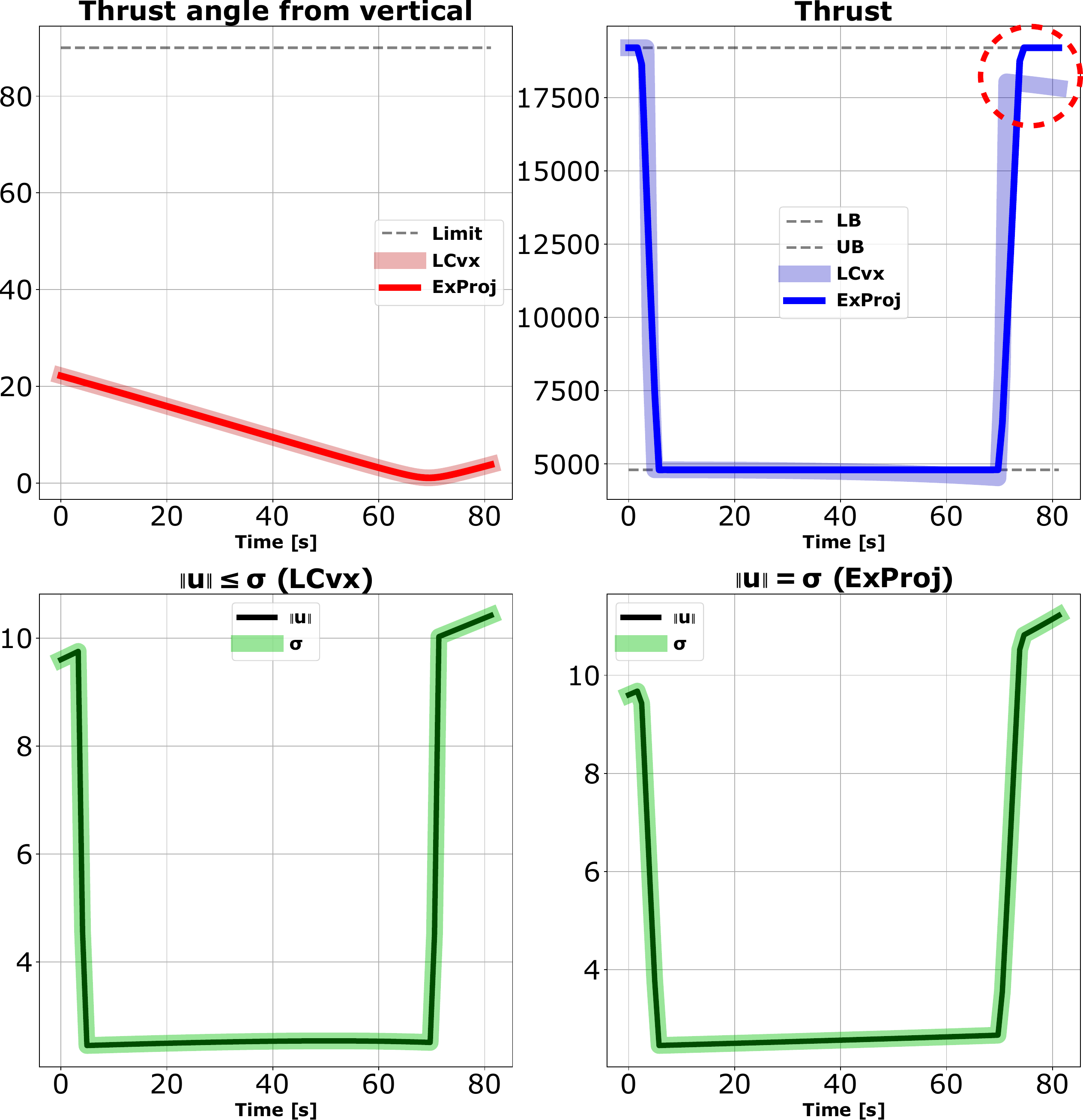}
\caption{\jhkhl{Optimization results for the scenario where the time of flight is strictly longer than the optimal time} ($t_f > t_f^*$). \jhkhl{Note that the effect of the linear approximation error on the thrust profile increases with the flight time (highlighted in dotted red circle where the LCvx solution fails to fully utilize the maximum thrust). This arises from the linear approximation for the maximum thrust constraint within the LCvx formulation.}} \label{fig:non_optimal_tof_more}
\end{figure}

For the second case, in which the time of flight is strictly longer than the optimal time ($t_f = 82\,\text{sec} > t_f^*$), the convexification is lossless with $\|u_i\|=\sigma_i$; however, the approximation error due to the linear approximation of the thrust bound constraints in \eqref{eqn:ncvx_thrust} means that the LCvx approach considers incorrect limits, and the discrepancy increases as the flight time increases.
This is clearly observed from the last part of the thrust profile in Figure~\ref{fig:non_optimal_tof_more}, where the LCvx solution uses significantly less thrust than the allowed maximum. Consequently, the LCvx solution is significantly conservative and thus suboptimal, whereas the ExProj solution is free of this suboptimality.

\begin{table}

\caption{Flight Test Parameters} \label{table:flight_test_parameters}
\begin{center}
\begin{tabular}{l c l l c}
\hline\hline
Parameter & Value & & Parameter & Value \\
\hline
$r_{\text{init}}$ (m) & (1.89, 4.59, 4.59) & & $m_{\text{wet}}$ (kg) & 0.03 \\
$\dot{r}_{\text{init}}$ (m/s) & (-0.18, 0.03, 0.03) & & $m_{\text{dry}}$ (kg) & 0.03 \\
$g$ (m/s$^2$)& (-9.81, 0, 0) & & $\rho_1$ (N) & 0.285\\
$\alpha$ (s/m) & $5\times10^{-9}$& & $\rho_2$ (N) & 0.306 \\
\hline
\hline
\end{tabular}
\end{center}
\end{table}

\section{Experimental Validation}
We briefly present the results obtained from flight tests conducted using a drone in the indoor flight arena equipped with motion capture systems. 

The proposed algorithm was implemented on an embedded GPU (NVIDIA Jetson AGX Orin) to efficiently compute and update the optimal PDG trajectory at a rate of 10 computations per second. Subsequently, the computed thrust command is transformed into both throttle and attitude commands of the drone.
In this test, we configured the parameter $\alpha$ in \eqref{eq:mass} to have a very small value, as the weight of the drone remains constant regardless of thrust usage. Furthermore, we incorporated a dynamic adjustment to the maximum allowable tilt angle of the thrust vector, which progressively decreases as the vehicle approaches the landing pad~\cite{adapt}. The parameters used for the flight test is given in Table~\ref{table:flight_test_parameters}.

\begin{figure*}
\begin{center}
    \includegraphics[width=\linewidth]{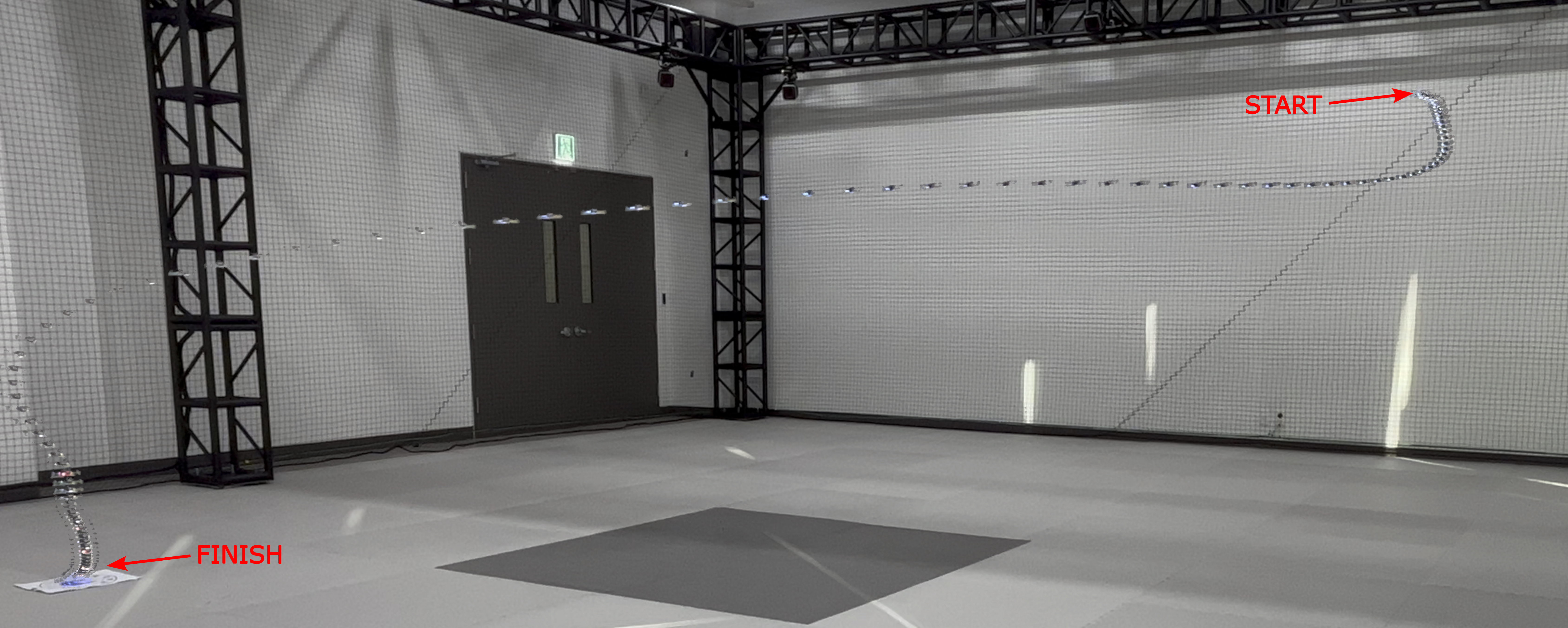}
\caption{Experimental validation of ExProj via indoor flight test. Sequence shot from the flight test is shown.} \label{fig:flight_test_real}
\end{center}
\end{figure*}

\begin{figure}
\begin{center}
	\includegraphics[width=0.88\linewidth]{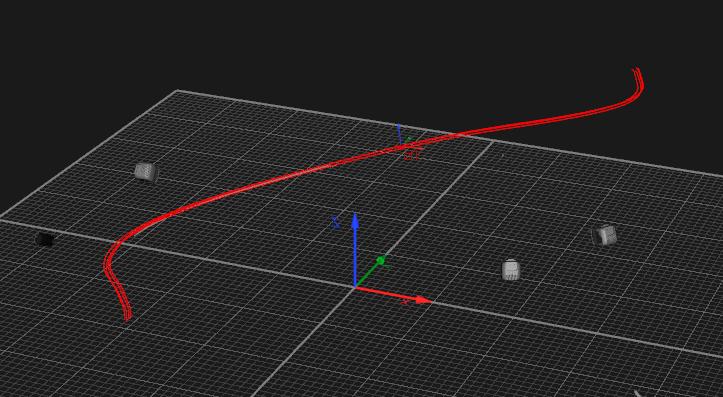}
	
	\footnotesize{(a) Trajectory obtained from the motion capture system.} 

	\bigskip

    \includegraphics[width=0.97\linewidth]{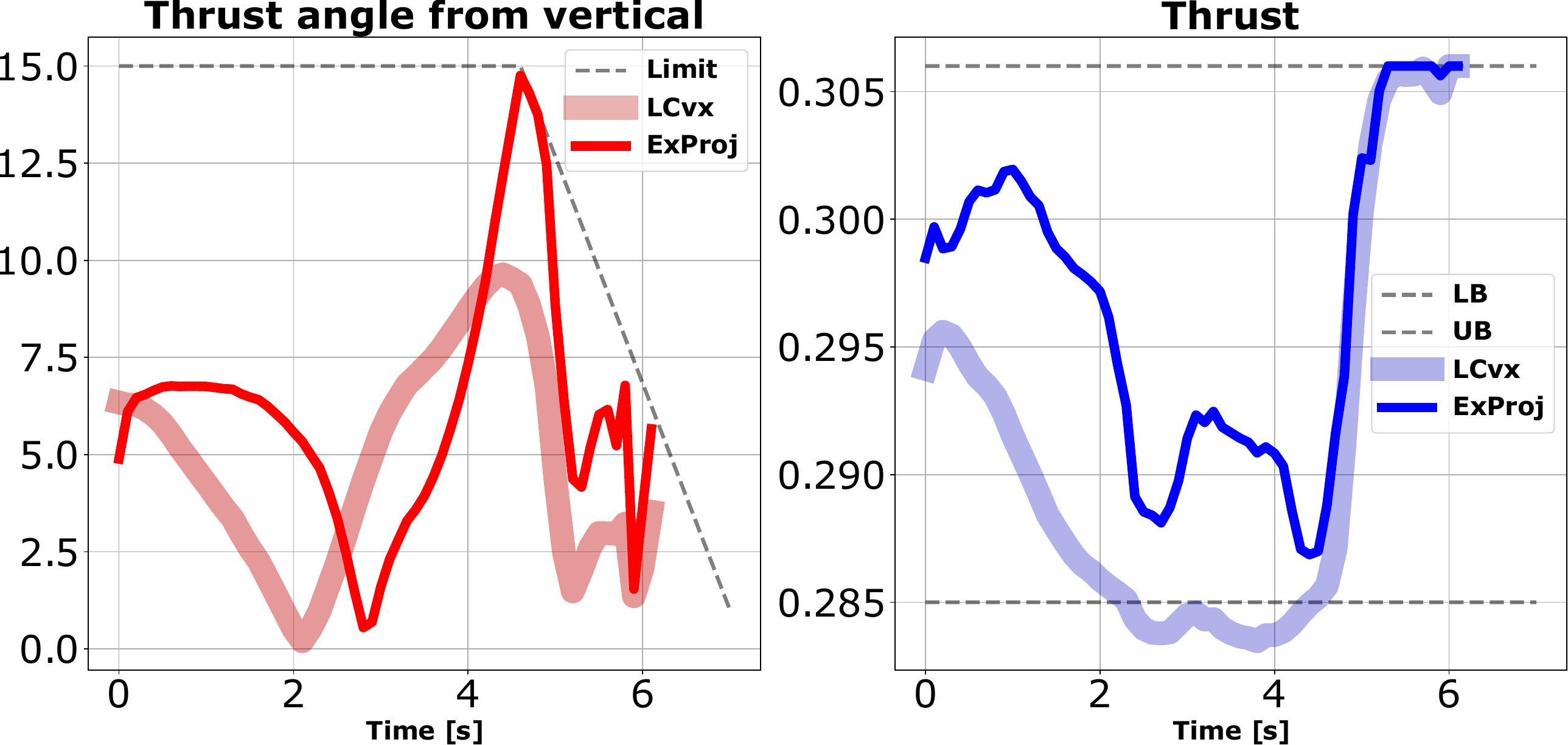}
    
	\footnotesize{(b) Tilt angle and the thrust profiles \jhkhl{from the indoor flight test}.} 

\caption{\jhkhl{Trajectory and the thrust profile from the indoor flight test. The flight results obtained through ExProj are compared to the computational results produced by LCvx, computed using the identical conditions. Note that the thrust profile derived from ExProj is feasible satisfying all the given thrust magnitude and angle bounds, while the thrust magnitude from LCvx solution violates the lower bound limit.}
}\label{fig:flight_test}
\end{center}
\end{figure}

In Figure~\ref{fig:flight_test_real}, we present the trajectory and a sequence shot obtained through the motion capture systems, and we present in 
Figure~\ref{fig:flight_test} the trajectory and thrust vector data acquired from the flight test, along with the corresponding tilt angle limit. The figure also includes the LCvx solution (depicted by lighter lines) computed from the identical initial conditions. Notably, it is observed that the LCvx solution is infeasible, necessitating a thrust level below the estabilished lower bound.

\section{Concluding Remarks}

In this paper, we proposed a first-order method that directly handles the nonconvexity arising in powered descent guidance problems.

Our approach combines first-order convex optimization algorithms with orthogonal projections onto nonconvex sets, which can be expansive.
Through a series of numerical examples, we verified the performance of the proposed algorithm and compared it with the most well-known standard convexification approach.
To the authors' best knowledge, it is the first approach that directly handles the nonconvex constraints in the PDG problem in the convex optimization frameworks and generates good feasible solutions in a variety of cases even when the existing standard approach fails.

\jhkhl{Numerical experiments reveals that our approach matches the standard convexification technique for lossless cases, outperforms in generating feasible solutions even when the standard techniques fail, and excels in terms of fuel consumption when the solution obtained from the standard approach is significantly suboptimal.} 

\jhkhl{Furthermore, we provided a concise overview of results obtained from an indoor flight test, demonstrating the efficacy of the proposed algorithm. It is observed that the solution derived from the proposed approach adeptly guided the vehicle with the desired flight time while maintaining the flight stability.}

\jhkhl{Although some results on the convergence of first-order methods on nonconvex problems have been reported recently, general results are not yet known. Therefore, a natural extension of the present work will be to analyze the convergence of this specific algorithm.}


\bibliographystyle{IEEEtran}
\bibliography{j9}
\end{document}